\documentclass[onecolumn, 12 pt, doublespace, fullpage, letterpaper]{article}

\usepackage{amsmath}
\usepackage{geometry}                		
\usepackage{graphicx}				
\usepackage{amssymb}
\usepackage{amsmath}	
\usepackage{amssymb}	
\usepackage{amsthm}	
\usepackage{geometry}
\usepackage{graphicx}
\usepackage{cite}

\usepackage{amssymb}
\usepackage{textcomp}
\usepackage{graphicx}
\usepackage{graphics}
\usepackage{epsfig}
\usepackage{epstopdf}
\usepackage{float}
\usepackage{color}
\definecolor{mygreen}{RGB}{28,172,0} 
\definecolor{myviloet}{RGB}{170,55,241}

\usepackage{dsfont} 

\usepackage{rotating}

\newtheorem{assumption}{Assumption}

\begin{document}

\title{An adaptive observer for wave equation's source estimation}
\author{Sharefa Asiri,$^1$ Taous-Meriem Laleg-Kirati,$^2$ and Chadia  Zayane-Aissa$^3$ \\
{\small $^{1,2,3}$Computer, Electrical and mathematical Sciences and Engineering, KAUST, Thuwal, Saudi Arabia.}\\
{\small Email addresses: $^1$sharefa.asiri@kaust.edu.sa}\\
{\small {\color{white}Email addresses:} $^2$taousmeriem.laleg@kaust.edu.sa}\\
{\small  {\color{white}Email addresses:} $^3$chadia.zayane@kaust.edu.sa}
}
\date{}
\maketitle

\section*{Abstract}
Observers are well known in control theory. Originally designed to estimate the hidden states of dynamical systems given some measurements, the observers scope has been recently extended to the estimation of some unknowns, for systems governed by partial differential equations.  In this paper, observers are used to solve inverse source problem for a one-dimensional wave equation. An adaptive observer is designed to estimate the state and source components for a fully discretized system. The effectiveness of the method is emphasized in noise-free and noisy cases and an insight on the impact of measurements' size and location is provided.

\section{Introduction} 
In this paper we are interested in an inverse source problem for the wave equation.  Inverse problems are usually solved using optimization techniques, where an appropriate cost function  is minimized. However the ill-posedness of such problems generates instability. Regularization techniques are then used to restore the stability. Among the regularization techniques, Tikhonov regularization \cite{GrCh:84} is probably the most used one. For instance, it has been applied to the wave equation in \cite{ChLi:83}  and  \cite{HoBo:05}. Other techniques have been also proposed to solve inverse problems; for example, in \cite{Ok:11}, a new minimization algorithm has been proposed to solve an inverse problem for the wave equation with unknown wave speed. Either regularized or not, most of the proposed methods end up with an optimization step which generally turns to be computationally  heavy, especially in the case of  large number of unknowns, and may require an extensive storage.

The objective of this paper is to  present an alternative method, based on observers, to solve the inverse source problem for the wave equation. Observers  are  well-known in control theory for state estimation in finite dimensional dynamical systems.  Presenting the distinctive feature and main advantage of operating recursively on direct problems, observers are gaining more and more interest in a wide variety of problems, including partial differential equations (PDE) systems. For instance, in \cite{MoChTa:08} states and parameters are estimated using an observer depending on a discretized space  for a mechanical system. In \cite{RaTuWe:09}, the initial state of a distributed parameter system has been estimated using two observers; one for the forward time and the other for the backward time.  A similar approach has been used in \cite{ChMi:10}, using the forward-backward approach to solve inverse source problem for the wave equation. An adaptive observer was applied in  \cite{GuoGuo:11} for parameter estimation and stabilization of one-dimensional wave equation where the boundary observation suffers from an unknown constant disturbance. A similar work was proposed  in \cite{GXHamm:12},with the state as  unknown and the boundary observation suffers from an arbitrary long time delay.

Dealing with PDEs, either with observers or classical inverse problems methods, poses the challenge of approximating infinite dimensional systems.  As regards observers, we can distinguish three approaches for studying such systems. The first approach considers the design of the observer in the continuous domain which requires a  semi-group theory analysis \cite{TuWe:09}. The second approach consists in the semi-discretization of the equation in space. The result of this semi-discretization can be usually written in the standard state-space representation in the continuous domain (in time) which makes the extension of the known methods easier. The third approach is the full-discretization of the PDE in space and time. In this case we can write the system in a discrete state-space representation and the methods known for such systems can be applied. We have chosen this latter approach as a straight forward application of observers.  We show that it can give good results provided that some conditions, aimed to minimize the effect of numerical errors resulting from discretization, are met.

 Another challenge, related to solving inverse problems in general, arises when it comes to measurement constraints. Indeed, from a practical point of view, we usually  do not have enough measurements to estimate all the unknowns. Dealing with this source of ill-posedness, means, in observers theory framework, satisfying the equivalent property of observability. Indeed, given the PDE system together with the measurements, we can test in a prior step  whether the unknown variables can be estimated fully or partially, regardless of the kind of observer to be used. For instance, in \cite{MoChTa:08,RaTuWe:09,GuoGuo:11,GXHamm:12}, the measurements were taken as the time derivative of the solution of  the wave equation. This kind of measurements gives a typical observability condition which has a positive effect on the stabilization, but  it is  less readily available than field measurements. Hence, some authors sought to solve inverse problems for wave equation using observers based on partial filed measurements, i.e. measurements taken from the solution of the wave equation, as in \cite{Ch:11}, \cite{ChCiMo:12}, and \cite{ChCiBuMo:12}.

In this paper, we consider a fully discretized version of a one dimensional wave equation and  we propose to apply the  adaptive
observer presented in \cite{GuZh:03} for the joint estimation of the states and the source term from partial measurements of the field. Adaptive observers are widely used in control theory for parameter estimation in adaptive control or fault estimation in fault detection and isolation \cite{PaCh:97}.
 In Section 2, problem statement is detailed. Then, the observer design is presented in Section 3. Finally,  numerical results are presented and discussed.

\section{Problem Statement}
Consider the one-dimensional wave equation with Dirichlet boundary conditions:
\begin{equation}
\left\{
\begin{array}{l}\label{wave equation}
 u_{tt}(x,t) - c^2u_{xx}(x,t)=f(x),\quad \quad x\in [0,l], \quad t\in [0,T], \\
   u(0,t)=0, \ \  u(l,t)=0,\\
 u(x,0)=r_1(x),\ \  u_t(x,0)=r_2(x),
\end{array}
\right.
\end{equation}
where $x$ is the space coordinate, $t$ is the time coordinate,
$r_1(x)$ and $r_2(x)$ are the initial conditions in  $\mathcal{L}^2 [0,l]$, $f(x) \in \mathcal{L}^2 [0,l]$ is the source function 
 which is assumed, for simplicity, to be independent of time. $c$ is the velocity which is known.
The notations $u_a$ and $u_{aa}$ refer to the first and second derivatives of $ u$ with respect to $a$, respectively. 

We first propose to rewrite  (\ref{wave equation}) in a system of first order PDEs by introducing two auxiliary variables
$ v(x,t)=u(x,t) $ and  $w(x,t)=u_t(x,t)$ 
and let 
\begin{equation}
\xi(x,t)=
\left[
\begin{array}{c c}
v(x,t), w(x,t)
\end{array}
\right]^T.
\end{equation}
 Then (\ref{wave equation})   can be written as follows,
\begin{equation}\label{wave equation2}
\left\{
\begin{array}{l}
\dfrac{\partial \xi (x,t)}{\partial t}=\mathcal{A} \xi(x,t)+F, \\
v(0,t)=0,\ v(l,t)=0, \\
v(x,0)=r_1(x),\ v_t(x,0)=r_2(x),\\
z=\mathcal{H} \xi(x,t),
\end{array}
\right.
\end{equation}
where the operator $\mathcal{A}$ is given by
$\mathcal{A}=
\left(
\begin{array}{cc}
 0 & I   \\
 c^2 \frac{\partial^2}{\partial x^2}&  0   
\end{array}
\right)
$ ,
 $F=
\left(
\begin{array}{c}
0 \\
f
\end{array}
\right)
$, $z$ is the output, and $\mathcal{H}$ is the observation operator such that
$\mathcal{H}=[\mathcal{H}_0 \quad 0]$ where $\mathcal{H}_0$ is a restriction operator on the measured domain.

%

\section*{Discretization of the problem}\label{discretized version}
 Discretizing system (\ref{wave equation2})  using  implicit Euler scheme in time and central finite difference discretization for the space gives the following discrete state-space representation:
\begin{eqnarray}\label{my state}
\left\{
\begin{array}{l}
  \xi^{j+1}=G\xi^j+B f^j+b, \\
 z^{j}=H\xi^j, \\
 f^{j+1}=f^{j}, \\
 j=1,2,\cdots,N_t ;
\end{array}
\right.
\end{eqnarray}
where  \[G=
\left(
\begin{array}{cc}
\Delta t E+ I & \Delta tI    \\
 E & I     \\
\end{array}
\right);
\] \[
 E=\frac{c^2\Delta t}{(\Delta x)^2}
 \begin{pmatrix}
  -2 &1 &  &   &  \\
 1 & -2 & \ddots &\ & \\
   & \ddots  & \ddots & 1 \\
  &  &1 & -2
 \end{pmatrix}
; B=
\left(
\begin{array}{c}
(\Delta t)^2 I \\
\Delta t I
\end{array}
\right); 
H= \left(
\begin{array}{cc}
H_m & {\bf 0}
\end{array}
\right);
\]
 $
H_m=
\left(
\begin{array}{ccc}
  0& \cdots  & 0  \\
  \vdots&I_m  &  \vdots  \\
 0& \cdots  & 0
\end{array}
\right); $  $I_m$ is the identity matrix of dimension $m$ where $m$ refers to the number of measurements, 
  and  $b$ is a term that includes the boundary conditions such that
 $b=\left( \frac{c^2 (\Delta t)^2}{(\Delta h)^2} v_1^j  \quad {\bf 0}_{1 \times (N_x-2)} \quad \frac{c^2 (\Delta t)^2}{(\Delta h)^2} v_{N_x}^j \quad \frac{c^2 (\Delta t)}{(\Delta h)^2} v_1^j \quad {\bf 0}_{1 \times (N_x-2)}  \quad \frac{c^2 (\Delta t)}{(\Delta h)^2} v_{N_x}^j \right)^T.$

This system is linear multiple-input multiple-output discrete time-invariant. If $N_x$ refers to the  space grid size, and $m$ refers to the number of measurements; thus, the state matrix $G$ is of dimension $2N_x \times 2N_x$,  the observer matrix $H$ is of dimension $m \times 2N_x$ , and the input matrix $B$ is of dimension $2N_x \times N_x$.

\section{Observer Design}\label{Observer Design}
 We propose to use an adaptive observer for the joint estimation of the states $v$ and $w$ and the source $f$. This observer has been proposed in \cite{GuZh:03}, and it has been developed for joint estimation of the state and the parameters. However,  we propose to generalize the idea behind this observer to estimate the input considering each spatial sample of the input as an independent parameter. The adaptive observer is given by the following system of equations, 
\begin{equation}\label{observer}
\left\{
\begin{array}{l}
\hat{z}^j= H\hat{\xi}^j,  \\
\Upsilon^{j+1}=(G-LH)\Upsilon^j+B,  \\
\hat{ f}^{j+1}=\hat{ f}^{j}+\sigma {\Upsilon^j}^T  H^T (z^j-\hat{z}^j), \\
\hat{\xi}^{j+1}=G\hat{\xi}^j+B\hat{ f^j}+b+
L(z^j-\hat{z}^j)+ \Upsilon^{j+1} (\hat{f}^{j+1}-\hat{f}^j),
\end{array}
\right.
\end{equation}
where $L$ is the observer gain matrix of dimension $2N_x \times m$, $\hat{\xi}^j$ and $\hat{ f}^{j}$ are the state and source estimates respectively, $\Upsilon^j$ is a matrix sequence obtained by linearly filtering $B$, and $\sigma$ is a scalar gain
 satisfying the following assumption as in \cite{GuZh:03}: 
\begin{assumption}\label{two assumption}
The scalar gain  $\sigma$ satisfies:
\begin{enumerate}
\item $\|\sqrt{\sigma} H\Upsilon^j  \|_2 \le 1$.
\item $\dfrac{1}{\kappa} \sum_{i=j}^{j+\kappa-1} { {\sigma \Upsilon ^i}^T  H^T  H \Upsilon ^i } \ge \beta I $
for some constant $\beta>0$, integer $\kappa >0$, and all $j$.

\end{enumerate}

\end{assumption}


\section{Numerical  results}

%

To test the performance of the observer, we generated a set of synthetic data using the following parameters:   $\Delta x=0.01$ and $l=2$, $\Delta t=0.01$ and $T=100$. Thus, $N_x=201$ and $N_t=10001$. The velocity is chosen to be $c^2=0.9$, and the source equal to $ f(x)=3\sin(5x)$. The method was implemented in Matlab and the tests were run for two main cases: noise-free and noisy data sets. In the noise-corrupted case, a white  Gaussian random noises with zero means were added to the states and to the measurements  with standard deviations $\sigma_{\xi}=0.007816$ and  $\sigma_z=0.01044$, respectively.
The gain matrix is selected to have fast and accurate convergence of the observer. However, standard pole placement  fails in this case, due to the size of the state matrix and the restricted number of measurements. Thus,  we took advantage of the particular structure of $G$ to design  the gain $L$. Indeed the matrix $G$  is  sparse  so we selected $L$ to be also a sparse matrix. The number of unknown entries is then reduced  and we identified them  such that  the eigenvalues of $(G - LH) $ are inside the unit circle. 

\begin{figure}[H]
\centering
\begin{tabular}{lcc}
 & {\bf Full Measurements} & {\bf Partial Measurements} \\
\begin{sideways}{\quad \quad \quad \quad {\bf Free-Noise}} \end{sideways}  &
\epsfig{file=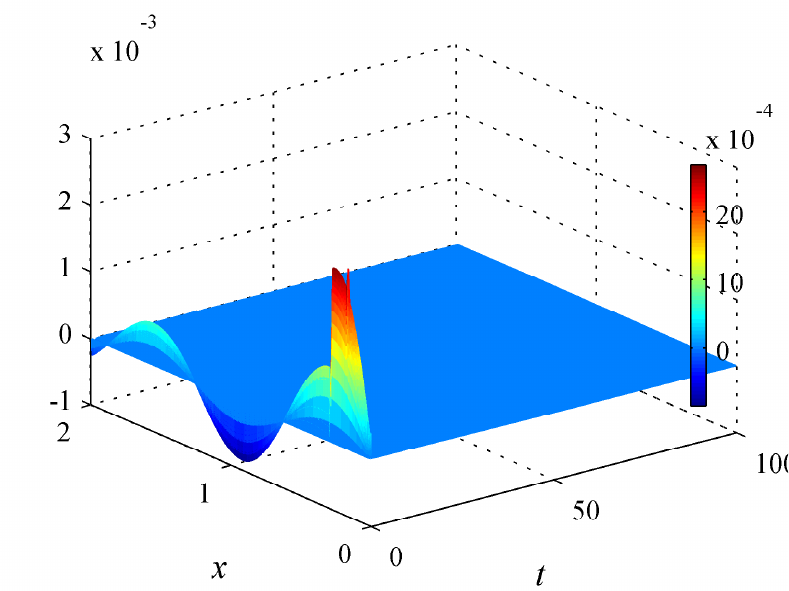,width=0.5\linewidth,clip=} &
\epsfig{file=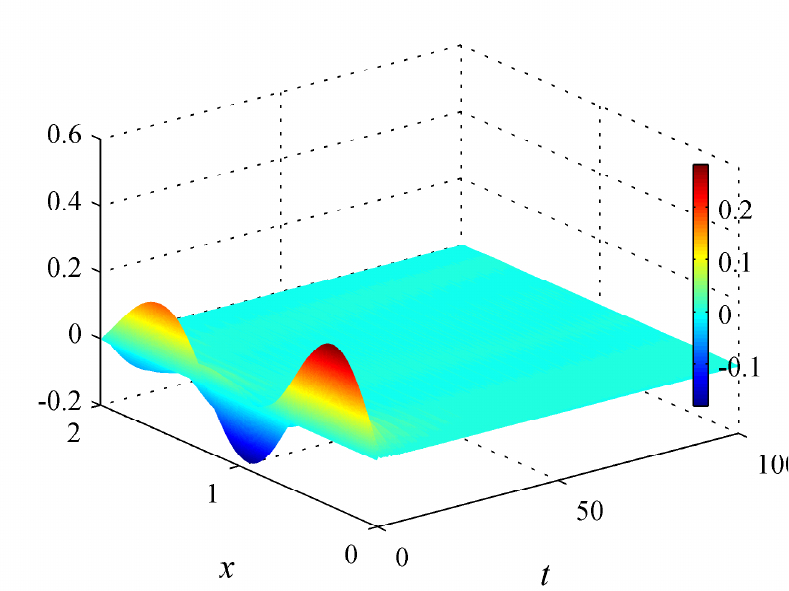,width=0.5\linewidth,clip=}\\
 & {\bf (a)} & {\bf (b)}\\
\begin{sideways}{\quad \quad \quad \quad {\bf Noise-Corrupted} } \end{sideways} &
\epsfig{file=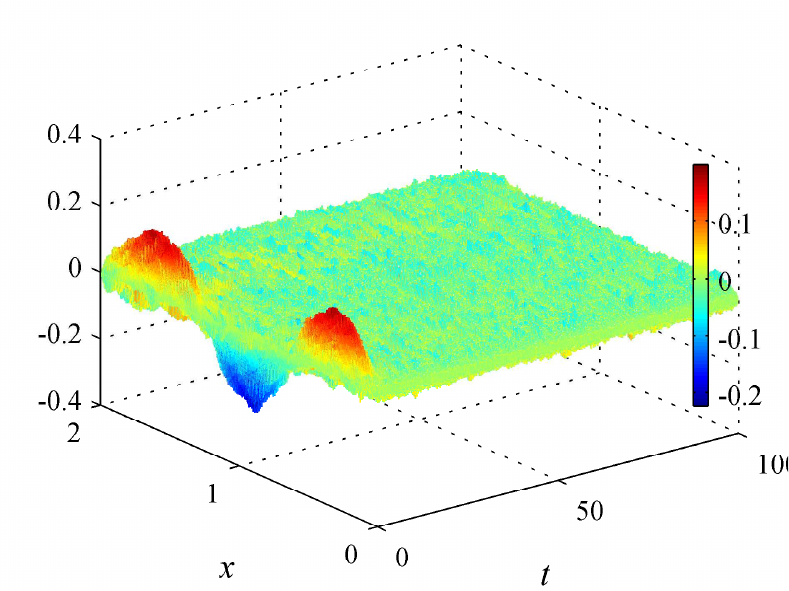,width=0.5\linewidth,clip=} &
\epsfig{file=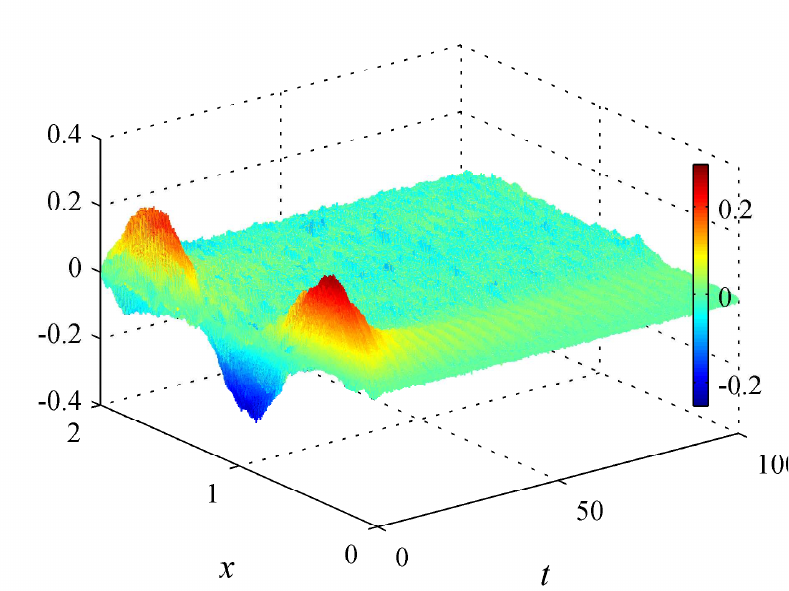,width=0.5\linewidth,height=0.38\linewidth,clip=}\\
 &{\bf (c)} & {\bf (d)}
\end{tabular}
\caption{The state error  ($\xi - \hat{\xi}$): (a) and (b) present the noise-free case with respect to full measurements and partial measurements, respectively. (c) and (d) show the noise-corrupted case with respect to full measurements and partial measurements, respectively. In the partial measurements cases, $50\%$ of the state components taken from the end.}
\label{state}
\end{figure}

\begin{figure}[H]
\centering
\begin{tabular}{lcc}
 & {\bf Full Measurements} & {\bf Partial Measurements} \\
\begin{sideways}{\quad \quad \quad \quad {\bf Free-Noise}} \end{sideways}  &
\epsfig{file=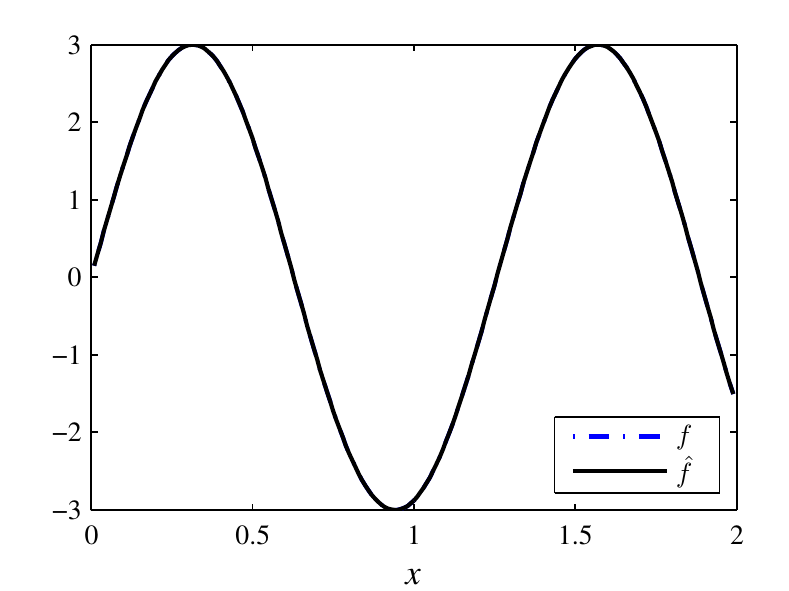,width=0.5\linewidth,clip=} &
\epsfig{file=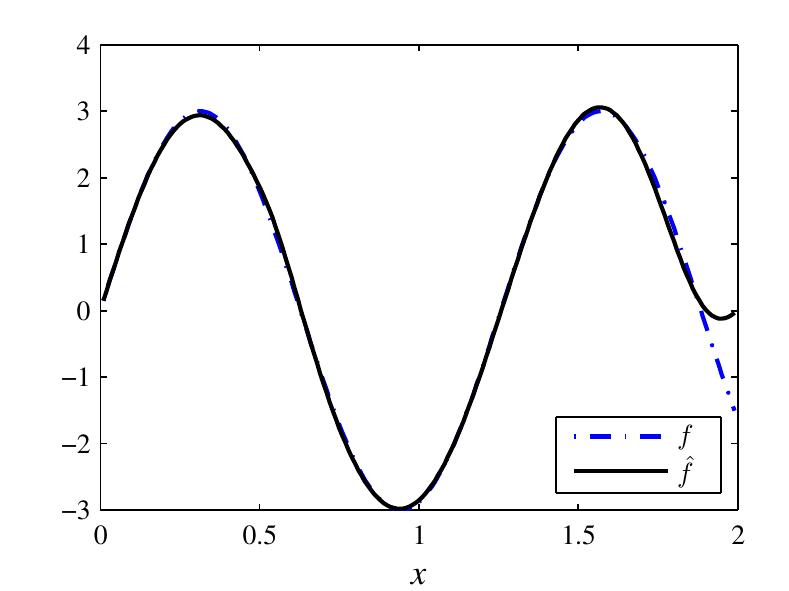,width=0.5\linewidth,clip=}\\
 & {\bf (a)} & {\bf (b)}\\
\begin{sideways}{\quad \quad \quad \quad {\bf Noise-Corrupted} } \end{sideways} &
\epsfig{file=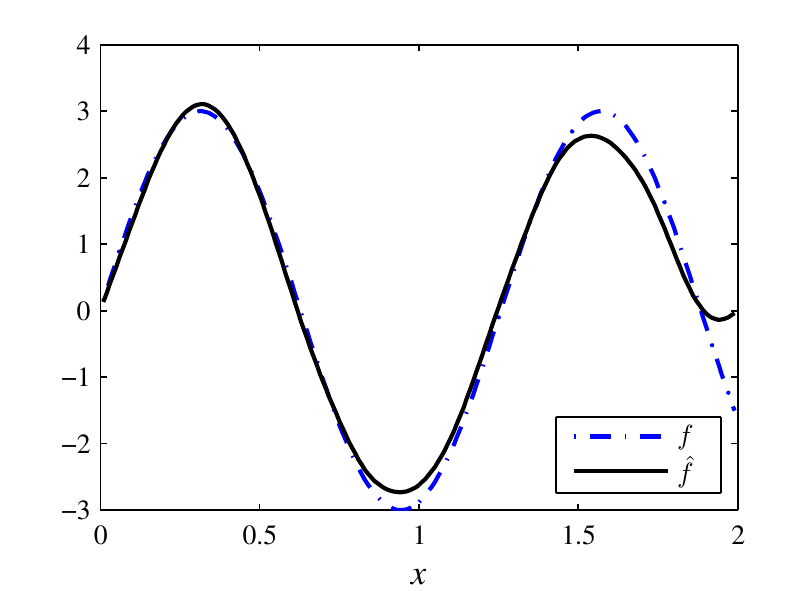,width=0.5\linewidth,clip=} &
\epsfig{file=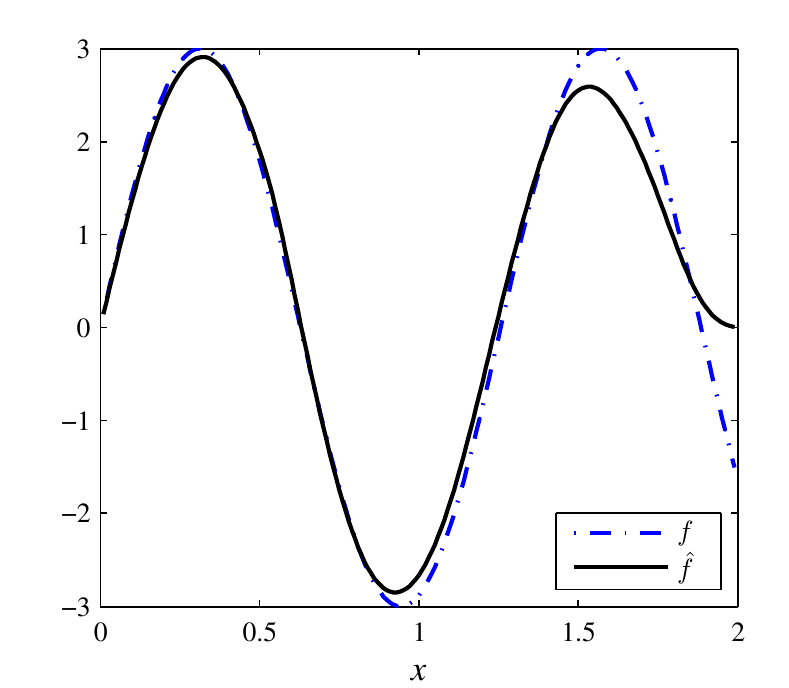,width=0.5\linewidth,height=0.38\linewidth,clip=}\\
 &{\bf (c)} & {\bf (d)}
\end{tabular}
\caption{The exact source $f$ (blue) and the estimated source $\hat{f}$ (black): (a) and (b) present the noise-free case with respect to full measurements and partial measurements, respectively. (c) and (d) show the noise-corrupted case with respect to full measurements and partial measurements, respectively. In the partial measurements cases, $50\%$ of the state components taken from the end.}
\label{source}
\end{figure}

Figure (\ref{state}) shows the error in the estimated state, and Figure (\ref{source}) presents the exact and the estimated source; both figures exhibit noise-free and noisy cases with respect to full and partial measurements. For the partial measurements, we supposed that the field is available on half of the space domain only. Table (\ref{source errors}) and Table (\ref{source errors_noisy}) show the minimum square error in the estimated source in noise-free and noisy cases, respectively. Both tables show the error in case of full measurements, partial measurements taken from the middle, and partial measurements taken from the end.

%
%
%


\begin{table}[H]
\caption{Source estimation errors in the noise-free case}
\begin{center}
\begin{tabular}{cc} \hline
Measurements  &MSE$=\sqrt{\dfrac{1}{N_x} \displaystyle \sum_{i=1}^{N_x} \left(f_i-\hat{f}_i \right)}$ \\ \hline
 Full  & $1.8168 \times10^{-14}$ \\ \hline
 Partial (middle)  &$ 0.3354$ \\ \hline
 Partial (end)  & $0.2096$\\ \hline
\end{tabular}
\end{center}
\label{source errors}
\end{table}%

\begin{table}[H]
\caption{Source estimation errors in the noisy case}
\begin{center}
\begin{tabular}{cc} \hline
Measurements  &MSE$=\sqrt{\dfrac{1}{N_x} \displaystyle \sum_{i=1}^{N_x} \left(f_i-\hat{f}_i \right)}$ \\ \hline
 Full  & $0.2865$ \\ \hline
 Partial (middle) & $0.4014$ \\ \hline
 Partial (end)  & $0.3213$\\ \hline
\end{tabular}
\end{center}
\label{source errors_noisy}
\end{table}

\subsection{Discussion}
In the noise-free case, the adaptive observer used in this paper provides a good estimate of the unknown source for the wave equation both  when the field is available on the whole space domain or when it is  available only on half the domain. However, we noticed in the second case a small error at the end of the interval.  In the noisy case also, the reconstruction is good but can be improved by a good choice of the gain $L$.   

We have studied the effect of number of measurements on the convergence of the proposed observer.   Obviously,  increasing number of measurements means increasing information about the state; thus, insuring the observability condition for all the states. However, for some applications, only few measurements can be available and the idea is to study the effect of this number on the convergence of the observer.

\begin{figure}[H]
\centering
\epsfig{file=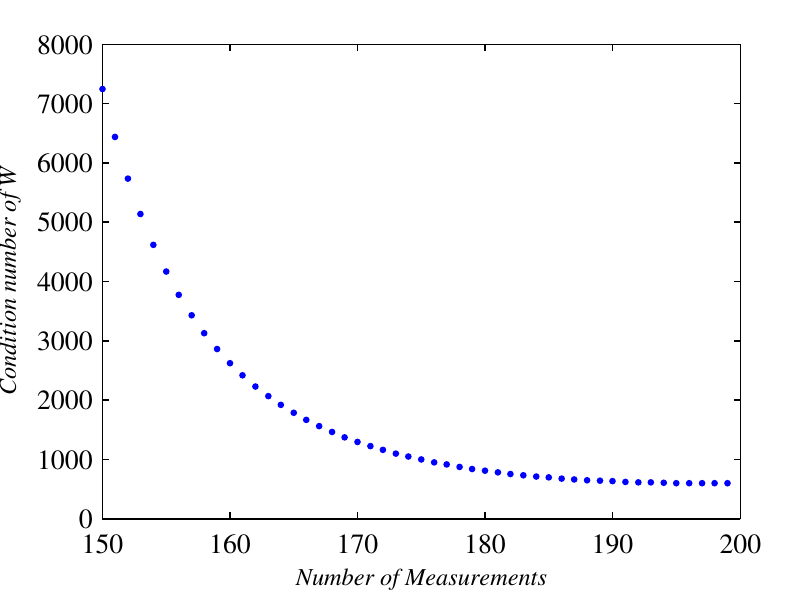,width=0.5\linewidth,clip=}
 \caption{Number of measurements versus the condition number of the observability matrix $W$.}\label{mesurements_vs_cond}
\end{figure}

The analysis of the error of estimation of the source with respect to the number of measurements shows that numerical issues may happen when we reduce the number of measurements below a threshold.  These numerical problems come in fact from the ill-conditioning of the observability matrix $W$, 
$W=\left(\begin{array}{ccccc}
 H &  HG & HG^2 & \cdots & HG^{n-1}
\end{array}
\right)^T.$ The decay of the condition number of the observability matrix $W$ as a function of the number of measurements is illustrated in Figure (\ref{mesurements_vs_cond}). It is well-known in control theory that the rank of $W$ gives the number of observable states.   It is known also that  a  high condition number for the observability matrix leads to a  nearly unobservable states \cite{Ch:91}.

\section{Conclusion}
In this paper, an adaptive observer for the joint estimation of the source and the states was designed. Numerical simulations for the source and states estimation using observer were presented, and they have proven the capability of observer to estimate both the source and the states.

\bibliographystyle{ieeetr}
\bibliography{References_3}

\end{document}